\documentclass{article}
\usepackage{graphicx}
\usepackage{color}
\usepackage{amsmath,amsfonts}

\begin{document}
\title{Positivity preserving high order schemes for angiogenesis models}
\author{A. Carpio (UCM), E. Cebri\'an (UBU)}

\maketitle


{\bf Abstract.} Hypoxy induced angiogenesis processes can be described coupling 
an integrodifferential kinetic equation of Fokker-Planck type with a diffusion equation
for the angiogenic factor. We propose high order positivity preserving schemes 
to approximate the marginal tip density by combining an asymptotic reduction
with weighted essentially non oscillatory  and strong stability preserving time 
discretization. We show that soliton-like solutions representing blood vessel formation 
and spread towards hypoxic regions are captured.

\section{Introduction}

\begin{figure}[h!]
\centering
\includegraphics[width=9cm]{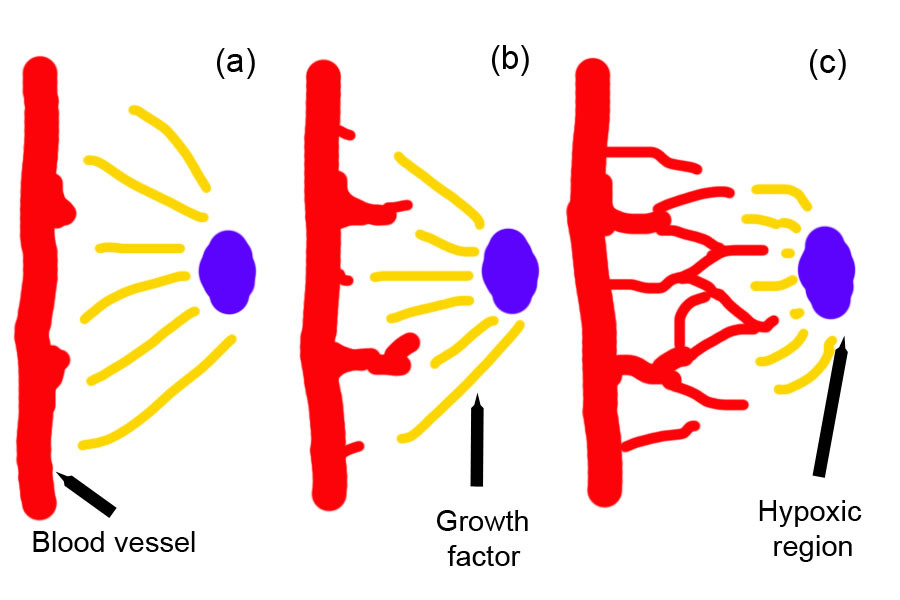}
\caption{Schematic illustration of the angiogenesis process: (a) The growth factor
released by the hypoxic region reaches an existing blood vessel. (b) Tip cells
lead the formation of new vessels. (b) The new vessels branch and spread towards
the hypoxic region.}
\label{fig1}
\end{figure}

Angiogenesis is a physiological process through which new blood vessels originate 
from pre-existing ones \cite{angioNat1,folkman}, through sprouting and splitting.
Blood vessel development is a normal and vital process in tissue growth, development 
and repair. However, it is also essential in the transition of tumors from a benign to a 
malignant state \cite{angioCh} and in many inflammatory and immune diseases.
Some treatments  rely on the use of angiogenesis inhibitors \cite{angioNat2,folkman,angioPr1,angioinhibitors}.
Angiogenesis is  related to hypoxia, that is, low levels of oxygen. Cells in hypoxic 
regions release sustances, such as VEGF (Vascular endothelial growth factor), which 
diffuse reaching adjacent blood vessels. VEGF  has been shown to increase the number 
of capillaries in a given network. In its presence, endothelial cells proliferate and migrate, 
eventually creating tubular structures resembling capillaries \cite{exercise}. Fig. \ref{fig1} 
illustrates the process.

Many models dealing with different features of angiogenesis processes have
been proposed, see \cite{angioComp1,angioComp2,angioComp3,angioPr2,angioComp4} 
for instance. Nevertheless, new experimental observations inspire novel mathematical 
models focusing on new features. The ultimate goal is to control the angiogenesis 
process to help healing as well as healthy tissue development.
Kinetic models allow us to gain insight on the dynamics of blood vessel networks 
by means of asymptotic and numerical tools. We consider a kinetic model which describes 
the stochastic nature of  blood vessel branching, derived from stochastic representations 
of the formation of blood vessel networks by means of ensemble averages  \cite{pre}. 
Let us denote by $p$ and $C$ the density of actively 
moving blood vessel tips  and the concentration of angiogenic factor released by hypoxic 
cells, respectively. Their  time evolution is governed by the following system of nondimensional equations:
\begin{eqnarray} \frac{\partial}{\partial t} p(\mathbf{x},\mathbf{v},t)\!&\!=\!&\!
 \alpha(C(\mathbf{x},t)) \delta_{\sigma_v}(\mathbf{v}\!-\!\mathbf{v}^0) p(\mathbf{x},\mathbf{v},t)  - \Gamma p(\mathbf{x},\mathbf{v},t) \int_0^t \!\! ds \!\! \int d{\bf v}' 
 p(\mathbf{x},\mathbf{v}',s)  \nonumber\\
 && - \mathbf{v}\cdot \nabla_\mathbf{x}   p(\mathbf{x},\mathbf{v},t) 
 + \beta {\rm div}_\mathbf{v} (\mathbf{v} p(\mathbf{x},\mathbf{v},t))+
\nonumber\\
&& - {\rm div}_\mathbf{v} \left[\beta\mathbf{F}\left(C(\mathbf{x},t)\right)p(\mathbf{x},\mathbf{v},t)  \right]\!  \frac{\beta}{2}
\Delta_\mathbf{v} p(\mathbf{x},\mathbf{v},t), \label{eq:p}  \\
\frac{\partial}{\partial t}C(\mathbf{x},t) &=& \kappa \Delta_{\mathbf x} C(\mathbf{x},t) - \chi C(\mathbf{x},t) j(\mathbf{x},t) \label{eq:c}, \\
p(\mathbf{x},\mathbf{v},0) &=& p_0(\mathbf{x},\mathbf{v}), \quad
C(\mathbf{x},0) =C_0(\mathbf{x}), \label{eq:pc0}
\end{eqnarray}
where
\begin{eqnarray} {
\alpha(C(\mathbf{x},t))=A\frac{C(\mathbf{x},t)}{1+C(\mathbf{x},t)}, \quad
 {\bf F}(C(\mathbf{x},t))= \frac{\delta_1}{(1+\Gamma_1C(\mathbf{x},t))^{q_1}}\nabla_{\mathbf x} C(\mathbf{x},t),}
\label{eq:alphaF} \\
j(\mathbf{x},t)= \int_{\mathbb R^N} 
\frac{|\mathbf{v}| }{ 1 + e^{(|\mathbf{v}- \mathbf{v}^0 |^2-\eta)/ \epsilon}}
p(\mathbf{x},\mathbf{v},t)\, d \mathbf{v},
\quad \rho(\mathbf{x},t)= \int_{\mathbb R^N} p(\mathbf{x},\mathbf{v},t)
\, d \mathbf{v}, \label{eq:intpintvp}
\end{eqnarray}
for ${\mathbf x} \in \Omega \subset  \mathbb{R}^N$, ${\mathbf v} \in  \mathbb{R}^N$, 
$N=2$, $t \in [0, \infty).$ The term $\alpha(C) \delta_{\sigma_v} p$ represents
the creation of new tips, while $- \Gamma p \int_0^t \rho ds$ quantifies tip destruction 
when existing vessels merge. 
The Fokker-Planck operator represents vessel branching and spread. The presence
of the force term $\mathbf F(C)$ in the equation (\ref{eq:p}) for the tip density $p$ 
forces spread towards the region where the source of growth factor is located. 
Similarly, the nonlocal sink $- C \, j$ in  (\ref{eq:c}) represents consumption of growth 
factor by the created vessel tips. 
The dimensionless parameters $\beta$, $\Gamma$, $\kappa$,  $\chi$, $A$, $\Gamma_1$, 
$\delta_1$, $\eta$, $\epsilon$ and $q_1$ are positive. $\delta_{\sigma_v}$ is a 
regularization of a delta function, a Gaussian centered about a reference velocity 
${\mathbf v_0}$:
\begin{equation}
\delta_{\sigma_v}(\mathbf{v}-\mathbf{v}_0)= \frac{1}{\pi\sigma_v^2}\, 
e^{-|\mathbf{v}-\mathbf{v}_0|^2/ \sigma_v^2}. \label{eq:deltav}
\end{equation}
An existence, stability and regularity theory in bounded and unbounded 
domains is developed in  \cite{bounded,unbounded}. Approximate soliton-like wave 
solutions have been found for the density of tips, which emerge at a pre-existing blood vessel and  advance until  the hypoxic region is reached \cite{pre16,jcp18}:
\begin{eqnarray} \begin{array}{l}
 \rho_s(x_1;K,c,X)=\frac{(2K\Gamma+\mu^2)c}{2\Gamma(c-F_1)}\mbox{sech}^2\!\left[\frac{\sqrt{2K\Gamma+\mu^2}}{2(c-F_1)}(x_1-X)\right]\!\!, \\
 \dot{X}=\frac{dX}{dt}=c,  \end{array} \label{soliton}
\end{eqnarray}
where  $c$  denotes the propagation speed and $X(t)$ the position of the maximum. 
Under some hypotheses, the parameters $K$ and $c$, satisfy a complementary
set of equations.

We  focus here in two dimensional settings, where these models may be adapted to describe angiogenesis disorders causing retinopathies \cite{retina}, for instance.  In three dimensions, the nonlocal terms must be modified to account for three dimensional peculiarities of blood vessel interaction and high order schemes would be necessary to reproduce such complex interactions. This task may involve a large computational  cost due to the presence of both space and velocity variables. Since the solutions of the model represent densities and concentrations, it is essential to use numerical discretizations ensuring positivity, a general issue with kinetic descriptions, which is challenging due to the presence of nonlocal terms and sources whose sign may change \cite{jcp18}. 

Finite difference and finite element schemes for Vlasov-Fokker-Planck models are 
considered in \cite{vpffe1,vpffe2,vfpfd}. Note that system (\ref{eq:p})-(\ref{eq:intpintvp}) 
resembles Vlasov-Poisson-Fokker-Planck (VPFP) equations except for the coupling 
with a diffusion problem and the presence of a nonlocal in time sink. In higher
dimensions, random particle methods \cite{vpfprandom} and deterministic particle 
methods combined with splitting of diffusion and convection operators \cite{vpfparticle}, 
have been proposed, as well as finite difference schemes combined with a change 
of variables \cite{vpfpfd}. These schemes are usually tested tracking convergence
to a stationary solution, under periodic boundary conditions to simplify.

Schemes for other kinetic problems exploit their specific features  \cite{russo}.  For 
Vlasov-Poisson (VP) and Vlasov-Maxwell (VM) systems, \cite{phd} gives an overview 
of possible approaches: WENO interpolation \cite{qiuweno}, discontinuous  Galerkin 
techniques \cite{sealdisgal}, conservative flux based methods \cite{filbertcons}, energy 
conserving finite differences \cite{filbertcomp}, particle in cell  techniques 
\cite{raviartpic,colellapic} and semi-lagrangian approaches \cite{semilagrangian} 
generalizing splitting methods \cite{knorrsplitting} are typical options. VM and VP 
models involve transport operators coupled to either Maxwell or Poisson problems.
They lack  the degenerate diffusion term part of Fokker-Planck (FP) operators.
In several dimensions, numerical schemes for Fokker-Planck equations are proposed 
in \cite{fpfd,fpc}, for instance; they often employ finite differences \cite{antonio}. 
Variants of Fokker-Planck operators are studied in \cite{degondfpl,taitanopfr}.

High order schemes for multidimensional VP and VPFP problems include  
approaches based on finite volume methods \cite{highfv},  
semi-lagrangian discretization \cite{highlag} or meshfree techniques \cite{highfree}.
In the slab geometry we consider for (\ref{eq:p})-(\ref{eq:intpintvp}) finite differences 
constitute a robust choice to preserve positivity and tackle the nonlocal in time 
coefficients, as well as the coupling with the diffusion equation. 
Our goal here is to devise high order schemes.
To lower the computational cost, we consider the limit of large $\beta$ and work
with an asymptotic approximation of the full model which provides reduced equations 
for key magnitudes: the marginal tip density and the angiogenic factor concentration.
We will focus on the final system,  which does not depend on the velocity,
combining positivity preserving weighted essentially nonoscillatory (WENO) 
discretizations of the differential operators in space and adequate strong stability 
preserving (SSP) discretizations in time \cite{shu,mstep,rk3,multistep}.

The paper is organized as follows. 
Section \ref{sec:2d} presents the two dimensional reduction of the full model. 
Sections \ref{sec:space} and  \ref{sec:time} detail the spatial discretization by means 
of WENO and SSP schemes. Section \ref{sec:num2d} displays numerical
solutions representing blood vessel tips migrating towards a hypoxic region. 
Finally, Section \ref{sec:conclusions} contains our conclusions.




\section{Asymptotic reduction}
\label{sec:2d}

in the limit of large friction, $\beta\to\infty$, the kinetic equation can be approximated by a reduced equation for the marginal density (\ref{eq:intpintvp}). The source terms in equation (\ref{eq:p}) (two first terms on its right-hand side) favor velocities in a small neighborhood of $\mathbf{v}^0$. Such velocities are those for which the tip creation term proportional to  $\alpha(C)\delta_{\sigma_v}(\mathbf{v} -\mathbf{v}^0)$ may balance the nonlocal tip destruction sink term $- \Gamma p \int_0^t \rho ds$. A reduced equation for $\rho$ is obtained by the Chapman-Enskog method \cite{pre16}. 

In the limit as $\beta\to\infty$, the marginal density $\rho(\mathbf{x},t)$
and the concentration $C(\mathbf{x},t)$ obey the equations \cite{jcp18}:
\begin{eqnarray}
\frac{\partial\rho}{\partial t}+{\rm div}_{\mathbf x}(\mathbf{F}\rho)-\frac{1}{2\beta}
\Delta_{\mathbf x} \rho=\mu\,\rho
-\Gamma\rho \int_0^t\rho(\mathbf{x},s)\, ds, \label{s1}\\
\mu=\frac{\alpha}{\pi}\left[1+\frac{\alpha}{2\pi\beta(1+\sigma_v^2)}\ln\!\left(1+\frac{1}{\sigma_v^2}\right)\!\right], \!\label{s2} \\ 
\frac{\partial}{\partial t}C(\mathbf{x},t)=\kappa \Delta_{\mathbf{x}} C(\mathbf{x},t)- \chi_1\, C(\mathbf{x},t)\,\rho(\mathbf{x},t),\label{s4}\\
\chi_1=\frac{\chi}{\pi}\int_0^\infty\int_{-\pi}^\pi \frac{\sqrt{1+V^2+2V\cos\varphi}}{1+e^{(V^2-\eta)/\epsilon}} e^{-V^2}V\,dVd\varphi. \label{s5}
\end{eqnarray}
To leading order, the marginal density and the density are related by
\begin{eqnarray}
p(\mathbf{x},\mathbf{v},t)\sim \frac{1}{\pi}e^{-|\mathbf{v}-\mathbf{v}_0|^2}\rho(\mathbf{x},t). \label{s3}
\end{eqnarray}
A positivity preserving order one scheme would follow by explicit forward time discretization, 
upwind treatment of transport terms, and centered schemes for Laplacians. Integral terms can be discretized using composite Simpson rules. 

To obtain a higher order scheme, we resort to positivity preserving WENO5 schemes for spatial operators, introduced in \cite{shu} for conservation laws with nonlinear fluxes and divergence free advection. Reference \cite{shu} starts from Lax-Friedrichs fluxes, we will use upwind fluxes here instead, and adapt the schemes to advection terms not necessarily divergence free. Also, the integral 
$I(\mathbf{x},t)\, = \int_0^t \rho(\mathbf{x},s)\, ds$ becomes an
additional equation 
\begin{eqnarray}  
I'(\mathbf{x},t) =  \rho(\mathbf{x},t), \quad  I(\mathbf{x},0)=0.  \label{dint}
\end{eqnarray}
The final system is composed of two  equations involving diffusion, reaction, and
convection terms, coupled to this additional ordinary differential equation.
In the next section, we detail the spatial discretization procedure.

\section{Weighted essentially non oscillatory space discretization}
\label{sec:space}

Let us consider a two dimensional convection-diffusion problem of the form
\begin{eqnarray}
u_{t} + (a u)_{x} + (b u)_{y}- d (u_{xx}+u_{yy})  = h(u), \label{basic}
\end{eqnarray}
where $d$ is a constant, either $d>0$ or $d=0$, and $\mathbf c=(a,b)$. 
In our case, $\mathbf c= \mathbf F= (F_1, F_2)$ depends on $\mathbf x = (x,y)$ 
and it is not necessarily divergence free.

Let us introduce a uniform rectangular mesh $(x_i,y_j)$, with
$x_{i+{1\over 2}} = x_{1 \over 2} + i \, \delta x$, $i=0, \ldots N_x$, and
$y_{j+{1\over 2}} = y_{1 \over 2} + j\, \delta y$, $j=0, \ldots N_y$. Let us
denote the midpoints by $x_i={x_{i-{1\over 2}} + x_{i+{1\over 2}}
\over 2}$ and $y_j={ y_{j-{1\over 2}} + y_{j+{1\over 2}} \over 2}$.
The spatial steps are given by $\delta x={X \over N_x},$ and
$\delta y= {Y \over N_y}$ for intervals $[0,X],$ $[0,Y].$ Several 
schemes for one dimensional conservations laws, and for
two dimensional nonlinear laws as well as advection equations
with divergence free non smooth velocities are described in \cite{shu}.
We specify here the procedure we use for our two dimensional
operators with non necessarily divergence free advection
coefficients, adapted from \cite{shu} with upwind fluxes, instead
of Lax-Friedrichs.

{\it Step 1.  Integrate (\ref{basic}) in the cell ${\cal C}_{ij}= [x_{j-{1\over 2}} x_{j+{1\over 2}}]
\times [y_{j-{1\over 2}} y_{j+{1\over 2}}]$.}
To handle the possible presence of a diffusive term, double cell averages
are used:
\begin{eqnarray}
\overline{\overline{u}}_{ij} = {1\over \delta x^2  \delta y^2}
\int_{y_{j-{1\over 2}}}^{y_{j+{1\over 2}}} \int_{y-{\delta y\over 2}}^{y+{\delta y\over 2}} 
\int_{x_{i-{1\over 2}}}^{x_{i+{1\over 2}}} \int_{x-{\delta x\over 2}}^{x+{\delta x\over 2}}
u(\xi,\eta) \, d \xi dx d\eta dy.  \label{double}
\end{eqnarray}
Integrating we find:
\begin{eqnarray} \begin{array}{ll}
{d \overline{\overline{u}}_{ij}  \over dt} =& - {1\over \delta x^2  \delta y^2}
\int_{y_{j-{1\over 2}}}^{y_{j+{1\over 2}}} \int_{y-{\delta y\over 2}}^{y+{\delta y\over 2}} 
\left[ \int_{x_i}^{x_{i+1}}  a u  \,dx  - \int_{x_{i-1}}^{x_{i}}  a u  \, dx  \right] d\eta dy 
\\
& - {1\over \delta x^2  \delta y^2}
\int_{x_{i-{1\over 2}}}^{x_{i+{1\over 2}}} \int_{x-{\delta x\over 2}}^{x+{\delta x\over 2}} 
\left[ \int_{y_j}^{y_{j+1}}  bu  \,dy  - \int_{y_{j-1}}^{y_{j}}  b u  \, dy  \right] d\xi dx 
\\
& + {d\over \delta x^2  \delta y^2}
\int_{y_{j-{1\over 2}}}^{y_{j+{1\over 2}}} \int_{y-{\delta y\over 2}}^{y+{\delta y\over 2}} 
[ u(x_{i+1},\eta) - 2 u(x_{i},\eta) + u(x_{i-1},\eta) ] \, d\eta dy
\\
&+ {d\over \delta x^2  \delta y^2}
\int_{x_{i-{1\over 2}}}^{x_{i+{1\over 2}}} \int_{x-{\delta x\over 2}}^{x+{\delta x\over 2}}
[ u(\xi,y_{j+1}) - 2 u(\xi,y_{j}) + u(\xi,y_{j-1}) ] \, d \xi dx
\\
& + {1 \over \delta x^2  \delta y^2} \overline{\overline{h(u)}}_{ij}.
\end{array} \label{basic_int} \end{eqnarray}

{\it Step 2: Replace the integrals by quadrature rules of $5$-th order.}
Three-point Legendre-Gauss quadrature rules are the standard choice. The weights 
and nodes are:
\begin{eqnarray}\begin{array}{llll}
x_\gamma= \left\{ -{\sqrt{15} \over 10}, 0, {\sqrt{15} \over 10} \right\}, &
w_\gamma= \left\{ {5\over 18}, {4\over 9}, {5 \over 18} \right\}, 
& \hbox{\rm in }   [-1/2,1/2], \\
x_{i+1}^\gamma= x_{i+1}+ x_\gamma \delta x, & w_\gamma,  
& \hbox{\rm in }   [x_i,x_{i+1}],
\end{array} \label{gaussl_2d} \end{eqnarray}
with $\gamma=1,2,3$. For these double integrals
\begin{eqnarray}\begin{array}{lll}
{1\over \delta x^2}
\int_{x_{i-{1\over 2}}}^{x_{i+{1\over 2}}} \int_{x-{\delta x\over 2}}^{x+{\delta x\over 2}}
v(\xi) d \xi dx  = 
{1\over \delta x} \sum_{\gamma=1}^3 \omega_\gamma \int_{x_i+ x_\gamma \delta x 
-{\delta x \over 2}}^{x_i+ x_\gamma \delta x +{\delta x \over 2}} u(\xi) d \xi  \\ =
 \sum_{\gamma=1}^3 \sum_{\gamma'=1}^3 \omega_\gamma \omega_{\gamma'}
v(x_i + x_\gamma \delta x + x_{\gamma'} \delta x)  = 
 \sum_{\gamma=1}^5 \tilde w_\gamma v(\tilde x_i^\gamma),
\end{array} \label{gaussl_2dq} \end{eqnarray}
where the final nodes $\tilde{x}_i^\gamma$ and weights $\tilde{w}_\gamma$ are
\begin{eqnarray}\begin{array}{llll}
\left\{ x_i - {\sqrt{15} \over 5} \delta x,  \; x_i - {\sqrt{15} \over 10} \delta x, \; 
x_i, \; x_i + {\sqrt{15} \over 10} \delta x, \; x_i + {\sqrt{15} \over 5} \delta x \right\},
\\
\left\{ w_1^2, \; 2 w_1 w_2, \; 2 w_1 w_3 + w_2^2, \; 2 w_3 w_2, \; w_3^2 \right\},
\end{array} \label{gaussl_2dwn} \end{eqnarray}
with $\gamma=1,\ldots,5.$ Using these quadrature rules, (\ref{basic_int}) becomes
\begin{eqnarray} \begin{array}{ll}
{d \overline{\overline{u}}_{ij}  \over dt} = & - {1\over \delta x }
\sum_{\alpha=1}^5 \sum_{\beta=1}^3 \tilde{w}_\alpha  {w}_\beta
[ a u(x_{i+{1\over 2}}^\beta,\tilde{y}_j^\alpha) - a u(x_{i-{1\over 2}}^\beta,\tilde{y}_j^\alpha) ]
 \\[1ex] 
& - {1\over \delta y }
\sum_{\alpha=1}^5 \sum_{\beta=1}^3 \tilde{w}_\alpha  {w}_\beta
[ b u(\tilde{x}_i^\alpha,y_{j+{1\over 2}}^\beta) - b u(\tilde{x}_i^\alpha,y_{j-{1\over 2}}^\beta) ]
 \\[1ex] 
& + {d\over \delta x^2} \sum_{\alpha=1}^5  \tilde{w}_\alpha
[ u(x_{i+1},\tilde{y}_j^\alpha) - 2 u(x_{i},\tilde{y}_j^\alpha) + u(x_{i-1},\tilde{y}_j^\alpha) ] 
 \\[1ex] 
& +  {d\over \delta y^2} \sum_{\alpha=1}^5  \tilde{w}_\alpha
[ u(\tilde{x}_i^\alpha,y_{j+1}) - 2 u(\tilde{x}_i^\alpha,y_{j}) + u(\tilde{x}_i^\alpha,y_{j-1}) ]
\\[1ex] 
& + \sum_{\alpha=1}^5 \sum_{\beta=1}^5  \tilde{w}_\alpha \tilde{w}_\beta h(\tilde{x}_i^\alpha,\tilde{y}_j^\beta).

\end{array} \label{basic_int2} \end{eqnarray}

{\it Step 3. Upwinding in the advection terms.} To ensure stability and preserve positivity,
it is convenient to discretize the advection terms in such a way that information
follows the characteristics \cite{shu,jcp18}. That is achieved  rewriting the terms
$a u(x_{i+{1\over 2}}^\beta,\tilde{y}_j^\alpha) $, $b u(\tilde{x}_j^\alpha,y_{i+{1\over 2}}^\beta)$ in terms of numerical fluxes that are consistent with the physical flux. In the absence
of diffusion and sources, this choice would also ensures that the scheme is conservative: speeds and possible shocks would be correctly captured. In this paper we use generalize upwind fluxes. Setting $f(u)= a u$, we  replace $a u(x_{i+{1\over 2}}^\beta,\tilde{y}_j^\alpha)$ by
\begin{eqnarray}
\hat{f}(u(x_{i+{1\over 2}}^\beta,\tilde{y}_j^\alpha)^-, u(x_{i+{1\over 2}}^\beta,\tilde{y}_j^\alpha)^+) = \left\{
\begin{array}{l} \hskip -2mm
a(x_{i+{1\over 2}}^\beta,\tilde{y}_j^\alpha) u(x_{i+{1\over 2}}^\beta,\tilde{y}_j^\alpha)^-, \;
a(x_{i+{1\over 2}}^\beta,\tilde{y}_j^\alpha) \geq 0, \\ \hskip -2mm
a(x_{i+{1\over 2}}^\beta,\tilde{y}_j^\alpha) u(x_{i+{1\over 2}}^\beta,\tilde{y}_j^\alpha)^+, \; a(x_{i+{1\over 2}}^\beta,\tilde{y}_j^\alpha) < 0, 
\end{array}\right.
\label{flux_upwind}
\end{eqnarray}
where $u(x_{i-{1\over 2}}^\beta,\cdot )^+$  and $u(x_{i+{1\over 2}}^\beta,\cdot )^-$ are the
approximations of $u(x_{i-{1\over 2}}^\beta,\cdot )$  and $u(x_{i+{1\over 2}}^\beta,\cdot )$ 
in cell ${\cal C}_{ij}$. Alternatively, we may use the Lax Friedrichs flux:
\begin{eqnarray}
\begin{array}{ll}
\hat{f}(u(x_{i+{1\over 2}}^\beta,\tilde{y}_j^\alpha)^-, u(x_{i+{1\over 2}}^\beta,\tilde{y}_j^\alpha)^+) \! = \!\! & {1 \over 2} a(x_{i+{1\over 2}}^\beta,\tilde{y}_j^\alpha) 
\left[u(x_{i+{1\over 2}}^\beta,\tilde{y}_j^\alpha)^- \!+\! u(x_{i+{1\over 2}}^\beta,\tilde{y}_j^\alpha)^+ \right]  \\[1ex] 
 & + {1\over 2} {\delta x \over \delta t} \left[u(x_{i+{1\over 2}}^\beta,\tilde{y}_j^\alpha)^-
- u(x_{i+{1\over 2}}^\beta,\tilde{y}_j^\alpha)^+\right].
\end{array} 
\label{flux_lf}
\end{eqnarray}
We denote by $\hat{g}$ the equivalent fluxes when $a$ is replaced by $b$ and
the roles of $x$ and $y$ are interchanged.
In this way we obtain the scheme
\begin{eqnarray} \begin{array}{ll}
{d \overline{\overline{u}}_{ij}  \over dt} = & - {1\over \delta x }
\left[ \hat{f}_{i+{1\over 2},j} - \hat{f}_{i-{1\over 2},j} \right] - {1\over \delta y }
\left[ \hat{g}_{i,j+{1\over 2}} - \hat{g}_{i,j-{1\over 2}}\right] \\[1ex]
& + {d\over \delta x^2} \sum_{\alpha=1}^5  \tilde{w}_\alpha
\left[ u(x_{i+1},\tilde{y}_j^\alpha) - 2 u(x_{i},\tilde{y}_j^\alpha) + u(x_{i-1},\tilde{y}_j^\alpha) 
\right] 
\\[1ex]
& +  {d\over \delta y^2} \sum_{\alpha=1}^5  \tilde{w}_\alpha
\left[ u(\tilde{x}_i^\alpha,y_{j+1}) - 2 u(\tilde{x}_i^\alpha,y_{j}) + u(\tilde{x}_i^\alpha,y_{j-1}) 
\right] 
\\[1ex]
& + \sum_{\alpha=1}^5 \sum_{\beta=1}^5  \tilde{w}_\alpha \tilde{w}_\beta h(\tilde{x}_i^\alpha,\tilde{y}_j^\beta), \end{array} \label{basic_int2} \end{eqnarray}
where
$\hat{f}_{i+{1\over 2},j}=\sum_{\alpha=1}^5 \sum_{\beta=1}^3 \tilde{w}_\alpha  {w}_\beta
\hat{f}(u(x_{i+{1\over 2}}^\beta,\tilde{y}_j^\alpha)^-, u(x_{i+{1\over 2}}^\beta,\tilde{y}_j^\alpha)^+)$
and
$\hat{g}_{i,j+{1\over 2}}=\sum_{\alpha=1}^5 \sum_{\beta=1}^3 \tilde{w}_\alpha  {w}_\beta
\hat{g}(u(\tilde{x}_i^\alpha,y_{j+{1\over 2}}^\beta)^-, u(\tilde{x}_i^\alpha,y_{j-{1\over 2}}^\beta)^+).$

{\it Step 4. Integrate in time.} The time discretization strategy for (\ref{basic_int2})
defines the order of the scheme in time. Forward Euler schemes are a simple
choice leading to order one, which provides stability and positivity. Known
the right hand $r_{ij}$ side of  (\ref{basic_int2}) at time $t_n$, the double
average $\overline{\overline{u}}_{ij}(t_{n+1}) = \overline{\overline{u}}_{ij}(t_{n})
+ \delta t \, r_{ij}(t_n)$. The next steps allow to approximate the values of $u$
at the nodes at time $t_{n+1}$ and the process is repeated at subsequent
times. Higher order time discretizations must belong to the SSP class (strong
stability preserving) as discussed in the next Section.

{\it Step 5}. The reconstruction of the required approximations of $u$ from the double averages  $\overline{\overline{u}}_{ij}$, is done dimension by dimension,
following a procedure detailed in\cite{shu}.  We recall the procedure in an Appendix
for ease of the reader.
Additionally, slope limiters can are enforced as described in \cite{shu}, if necessary.
The scheme is expected to be stable when
\begin{eqnarray}
({dt \over dx} + {dt \over dy} ) {\rm max}{|a|,|b|} \leq {1 \over 2} \hat{w}_1 {\rm min}
\{w_1,w_2,w_3\} = {1 \over 2} {1 \over 12} {5 \over 18} =1 \times 10^{-2}, \\
({dt \over dx^2} + {dt \over dy^2} ) {\rm max}{|d|} \leq {1 \over 4} \tilde{w}_3
= {1 \over 4} {114 \over 324} = 8 \times 10^{-2}.
\end{eqnarray}

\section{Strong stability preserving time discretization}
\label{sec:time}

To maintain positivity and stability, we work with strong stability  preserving 
(SSP) time discretizations. 
High order time discretizations of Runge-Kutta type increase the cost since 
several evaluations of the right hand side are required per time step. Multistep 
schemes avoid this problem reusing previous step, but increase storage needs 
and usually require smaller time steps to ensure stability. 
In spite of their theoretical order, WENO5 schemes might degenerate to order 
two in practice in smooth regions. Therefore, we will use SSP schemes or order 
2 or 3. Usual choices for third order accuracy are a third order SSP multistep  
method \cite{shu}
\begin{eqnarray} \begin{array}{l}
u(t_{n+1})= {16 \over 27} \left( u(t_n) + 3 \delta t  \, r(u(t_n)) \right)
+ {11 \over 27} \left( u(t_{n-3}) + {12 \over 11} \delta t \, r(u(t_{n-3}))
\right), 
\end{array} \label{mstep} \end{eqnarray}
and a third order Runge Kutta method \cite{rk3}
\begin{eqnarray} \begin{array}{l}
u^{(1)}= u(t_n) + \delta t \, r(u(t_n)), \\
u^{(2)}= {3 \over 2} u(t_n) + {1\over 4}  u^{(1)} 
+ {1\over 4} \delta t  \, r(u^{(1)}), \\
u(t_{n+1})= {1 \over 3} u(t_n) + {2\over 3}  u^{(2)} 
+ {2\over 3} \delta t  \, r(u^{(2)}). 
\end{array} \label{rk3} \end{eqnarray}
The stability of SSP discretizations is governed by a CFL number $c$ as
follows. If  Euler forward time discretization applied to equation $u_t = r(u)$ 
remains stable under the condition $\delta t \leq \delta t_0$, then 
SSP time discretization is stable when $\delta t \leq d \, \delta t_0$. For 
the multistep method  $d=1/3$ while $d=1$ ($d_{eff}=1/3$) for the RK3 (\ref{rk3}). 
For second order approximations, the RK2 scheme is
\begin{eqnarray}  \begin{array}{l}
u^{(1)}= u(t_n) + \delta t \, r(u(t_n)), \\
u(t_{n+1})= {1 \over 2} u(t_n) + {1\over 2}  u^{(1)} 
+ {1\over 2} \delta t  \, r(u^{(1)}), 
\end{array} \label{rk2} \end{eqnarray}
with $d=1$ ($d_{eff}=1/2$).
In our case, the spatial operator $r(u(t_n))$ is the operator obtained discretizing
the space variables, time excluded.   

The final scheme for our system would read as follows. 
\begin{itemize}
\item Initialize the scheme applying the RK3 scheme (\ref{rk3}) to the spatially
discretized versions of equations (\ref{s1}) and (\ref{s4}), according to Section
(\ref{space}), and to equation (\ref{dint}). Store the first three steps, for
$n=1,2,3$.
\item At each $n>3$, apply the multistep scheme (\ref{mstep}) to the spatially
discretized versions of equations (\ref{s1}) and (\ref{s4}), according to Section
(\ref{space}), and to equation (\ref{dint}). 
\end{itemize}

\section{Numerical tests in two dimensions}
\label{sec:num2d}

In this section, we present numerical solutions of the reduced equations
(\ref{s1})-(\ref{dint}) for appropriate values of the parameters as listed in Table 
\ref{table1}. The slab geometry of the advancing angiogenic network is schematized 
in Figure \ref{fig1}. More precisely, we consider a slab $(0,1)\times \mathbb{R}$, 
and set ${\mathbf x}= (x_1,x_2)$. To ease the implementation of WENO
discretizations, we slightly modify the boundary conditions proposed in the original 
model \cite{jcp18}.
For numerical purposes we truncate the spatial domain and set 
$\Omega_{\mathbf x}=[0,1] \times [-1.5,1.5].$ We impose zero Dirichlet 
boundary conditions  for $C$ at $x_2=\pm 1.5$ and $x_1=0$,
as well as zero boundary conditions for $\rho$.  At $x_1=1$, we use the
function
\begin{equation}
g(t,x_2) = c_L(t) e^{-a^2 x_2^2}, \quad 
t>0, x_2 \in \mathbb{R}, \label{bc:c}
\end{equation}
to generate boundary values. Here, $c_L(t)>0$ to represents the influx of 
angiogenic factor produced by the core of the hypoxic region and $1/a$ is a 
characteristic length thereof. This function will decrease with time as new blood
vessels reach the tumor.  We set in our simulations $a=1/0.3$ and $c_L(t)=1.1$.

\begin{table}[ht]
\begin{center}\begin{tabular}{cccccccccc}
 \hline
$\delta_1$ & $\beta$ &$A$& $\Gamma$& $\Gamma_1,q_1$ &$\kappa$&$\chi$&$\eta$
&$\epsilon$ & $\sigma_v$\\
0.255 & 5.88 &  22.42  & 0.135 & 1 &  0.0045  & 0.002 & 15 &  0.001 & 0.08 \\
 \hline
\end{tabular}
\end{center}
\caption{Dimensionless parameters. } 
\label{table1}
\end{table} 

Fot the initial data, we consider
\begin{eqnarray}
p(x_1,x_2,v_1,v_2,0)={2\over \pi^2} {1 \over 0.0048}
e^{-\left(x_1\over 0.06\right)^2} \sum_{j=1}^{20}
e^{-\left({x_2-x_2^j\over 0.08} \right)^2} 
e^{-[(v_1-v_1^0)^2+(v_2-v_2^0)^2]},  \label{ic:p}\\
C(x_1,x_2,0)=1.1 \, e^{-[({ x_1-1\over 1.5})^2 + ({ x_2\over 0.3})^2]}, \label{ic:c}
\end{eqnarray}
where the vessel tip locations $x_2^j=-0.3+ (j-1) {0.6 \over 19}$, $j=1,2,\ldots,20,$ 
are $20$ equispaced values in $[-0.3,0.3]$, $v_1^0= \cos(\pi/10)$ and 
$v_2^0=\sin(\pi/10)$. 
Initial values for $\rho(\mathbf x,0)$ are obtained integrating (\ref{ic:p}) over the 
truncated computational velocity domain $\Omega_{\mathbf v}= [- 4, 4]^2$.
The steps are $\delta x= 0.02$, 
and $\delta t = 2.2906 \times 10^{-4}$.

\begin{figure}[h!]
\centering
(a) \hskip 3cm (b) \hskip 3cm (c) \\
\includegraphics[width=4cm]{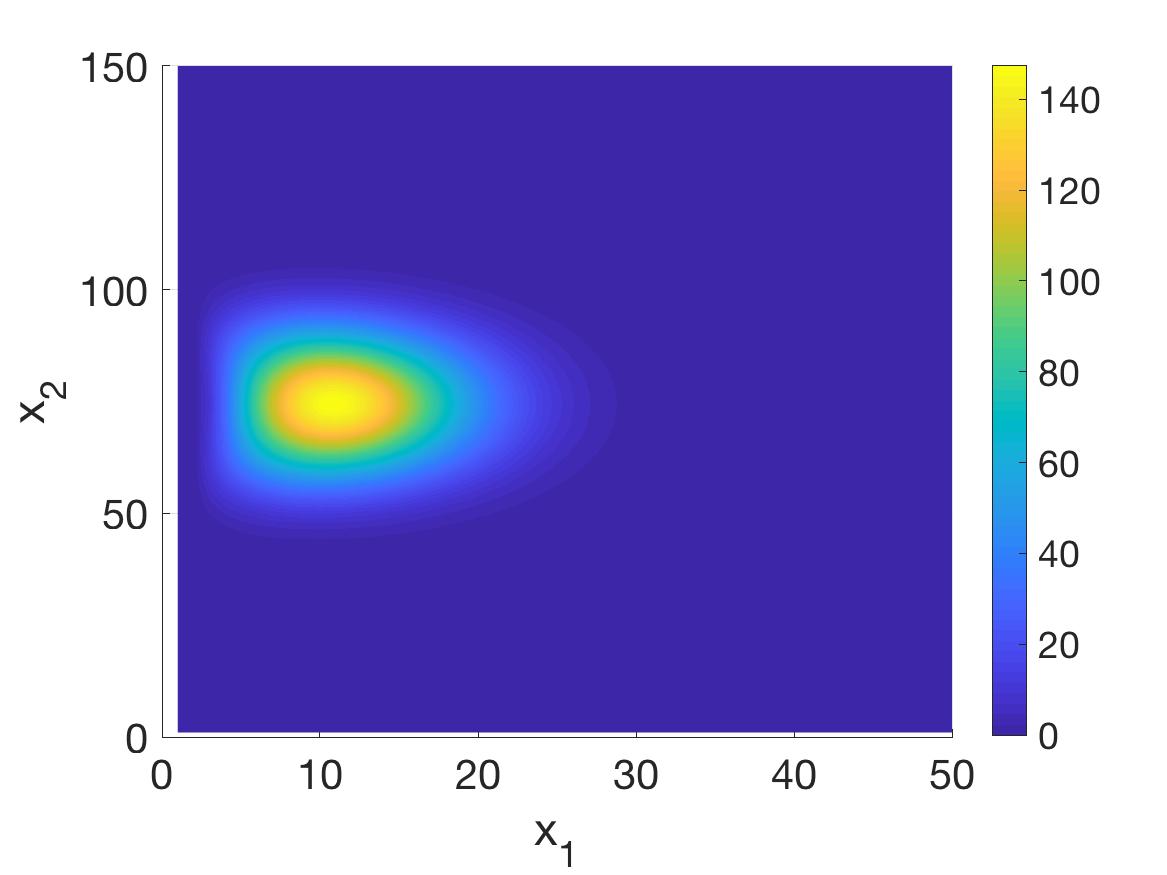} \hskip -5mm
\includegraphics[width=4cm]{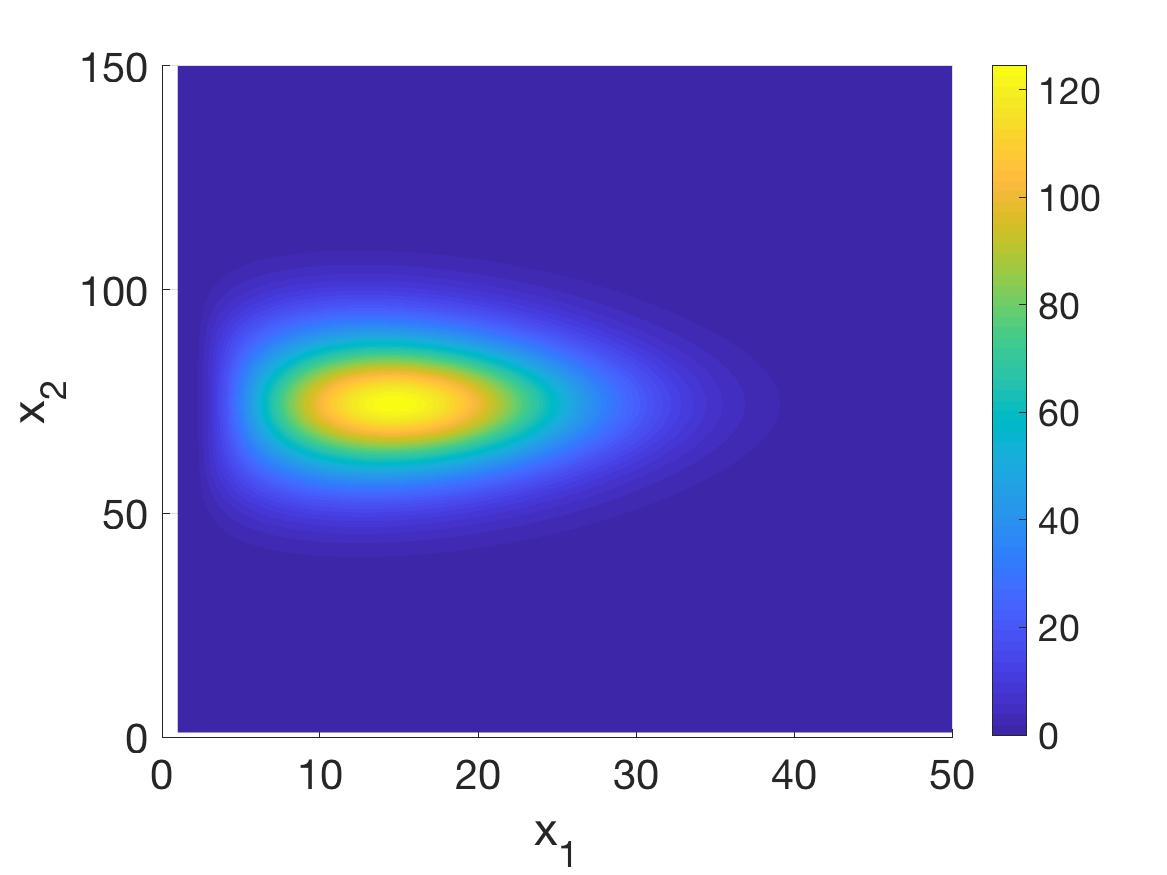} \hskip -5mm
\includegraphics[width=4cm]{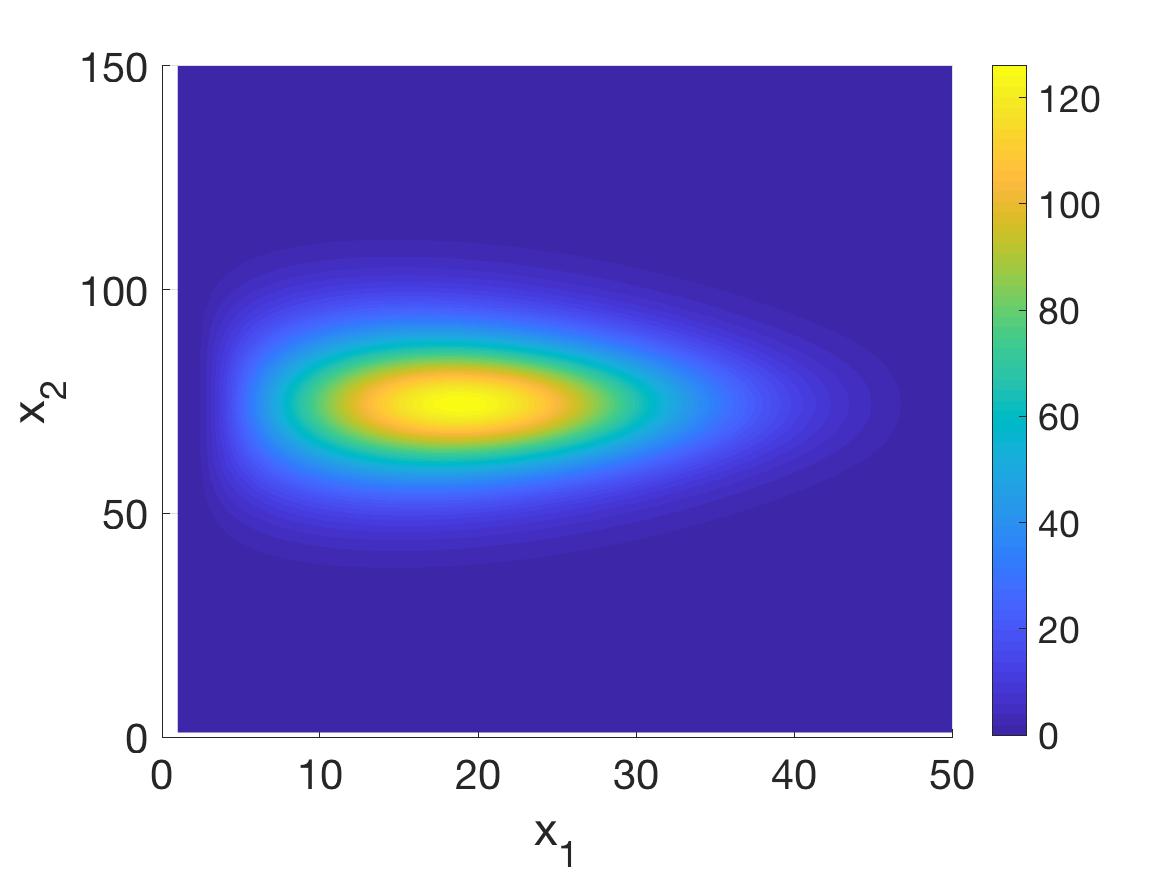} \\
(d) \hskip 3cm (e) \hskip 3cm (f) \\
\includegraphics[width=4cm]{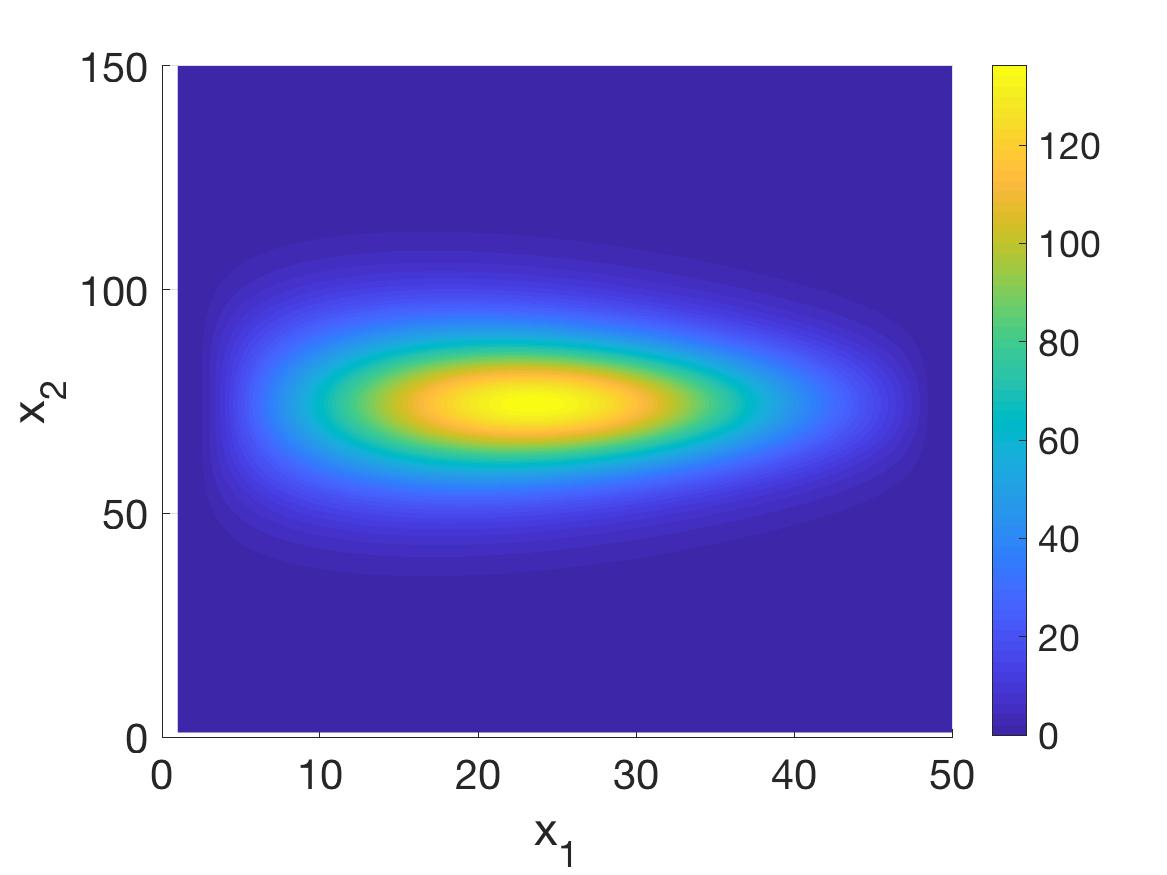} \hskip -5mm
\includegraphics[width=4cm]{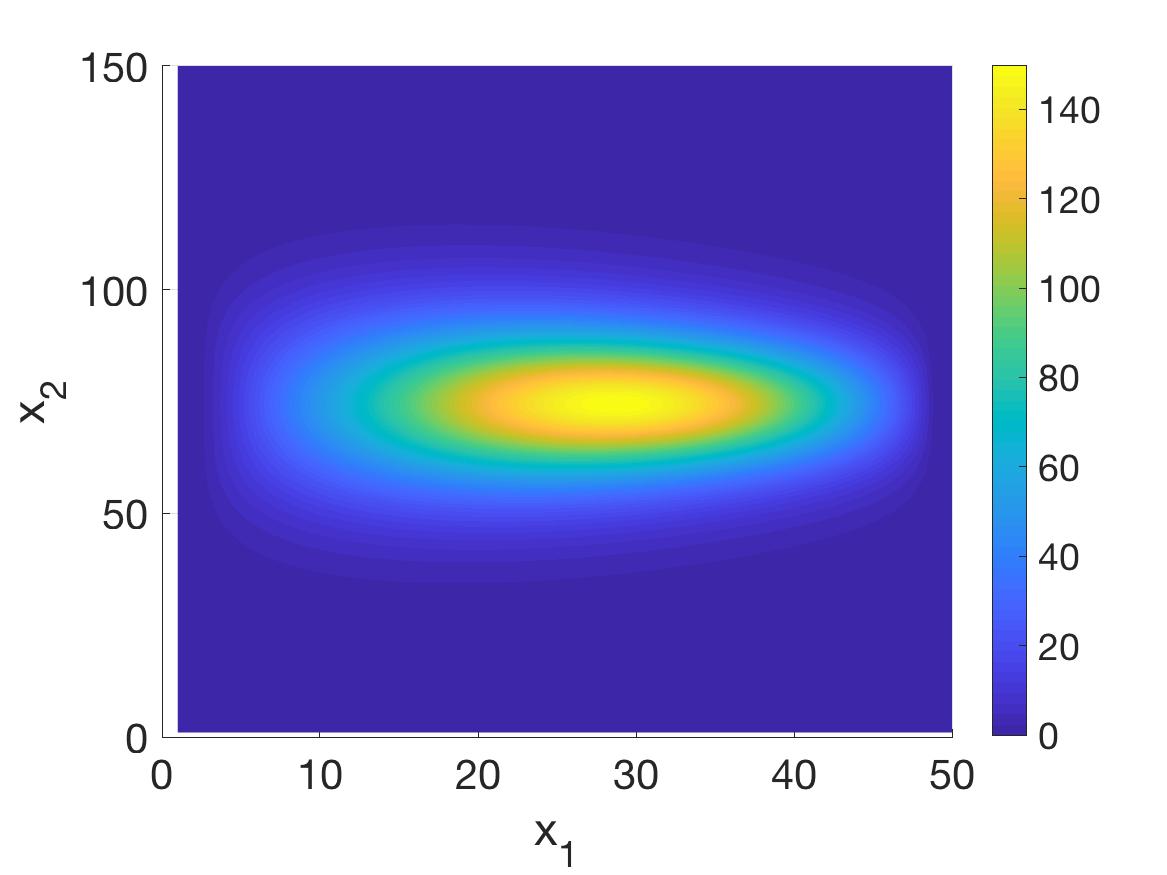} \hskip -5mm
\includegraphics[width=4cm]{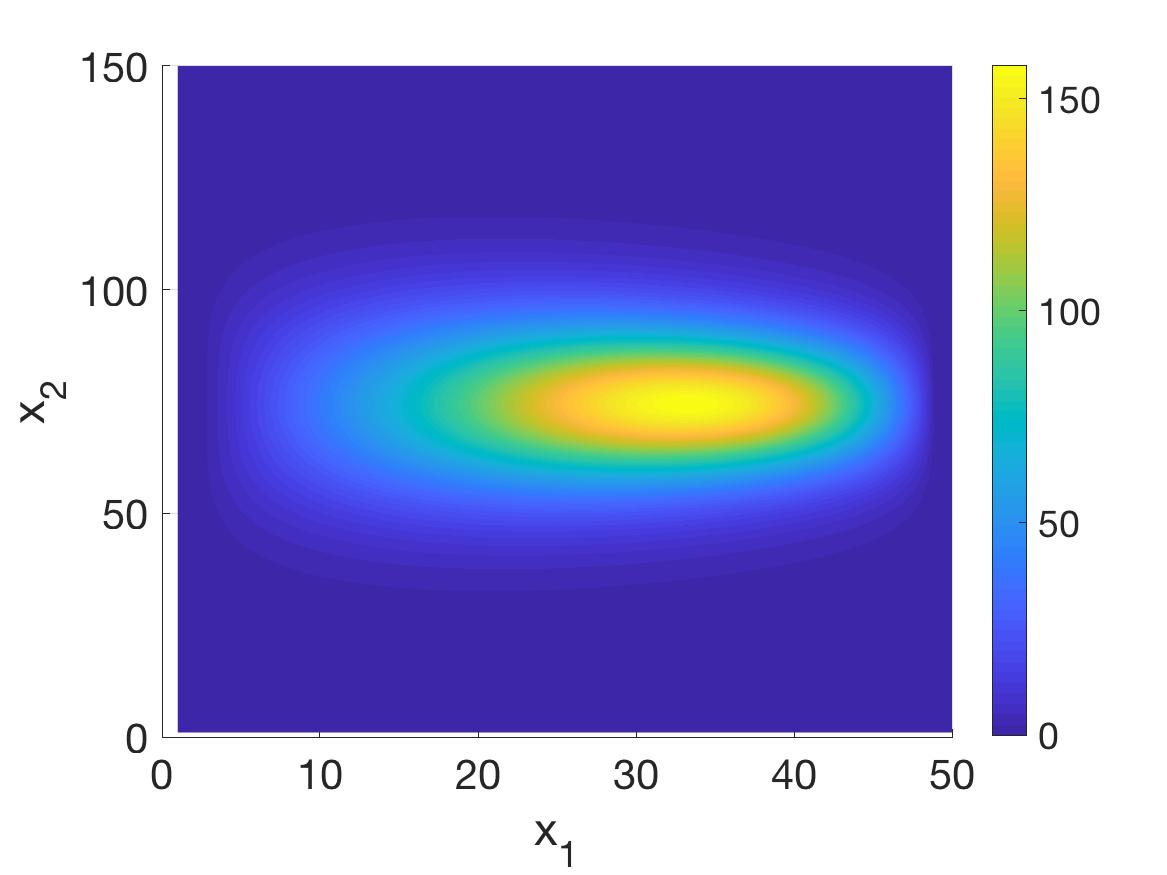} 
\caption{Marginal tip density obtained solving the two dimensional
reduced model at times
(a) $0.1145$, (b) $0.2290$, (c)  $0.3435$, (d) $0.4580$,
(e)  $0.5725$, (f) $0.6870$.
A soliton-like pattern of new tips forms which travels from the pre-existing 
vessel (left) to the hypoxic region releasing the angiogenic factor (right).} 
\label{fig2}
\end{figure}

\begin{figure}[h!]
\centering
(a) \hskip 3cm (b) \hskip 3cm (c) \\
\includegraphics[width=4cm]{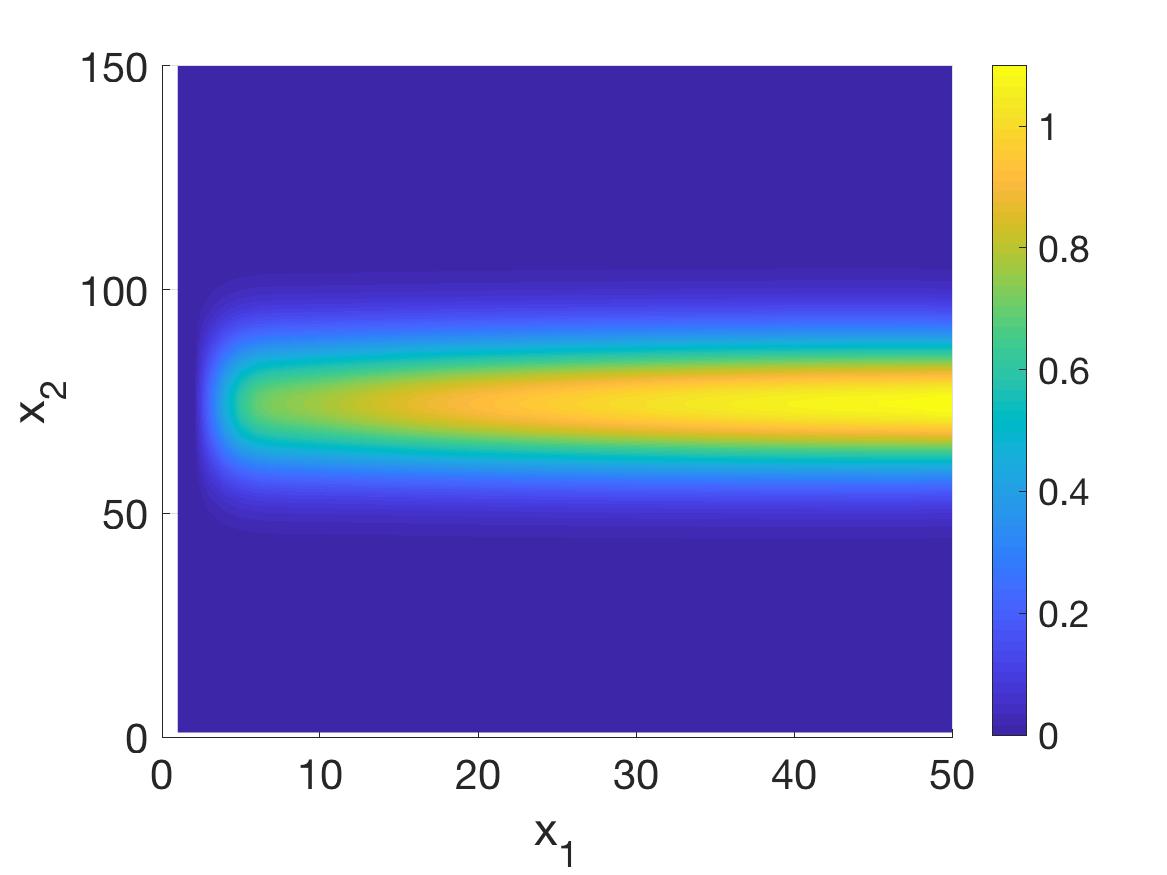} \hskip -5mm
\includegraphics[width=4cm]{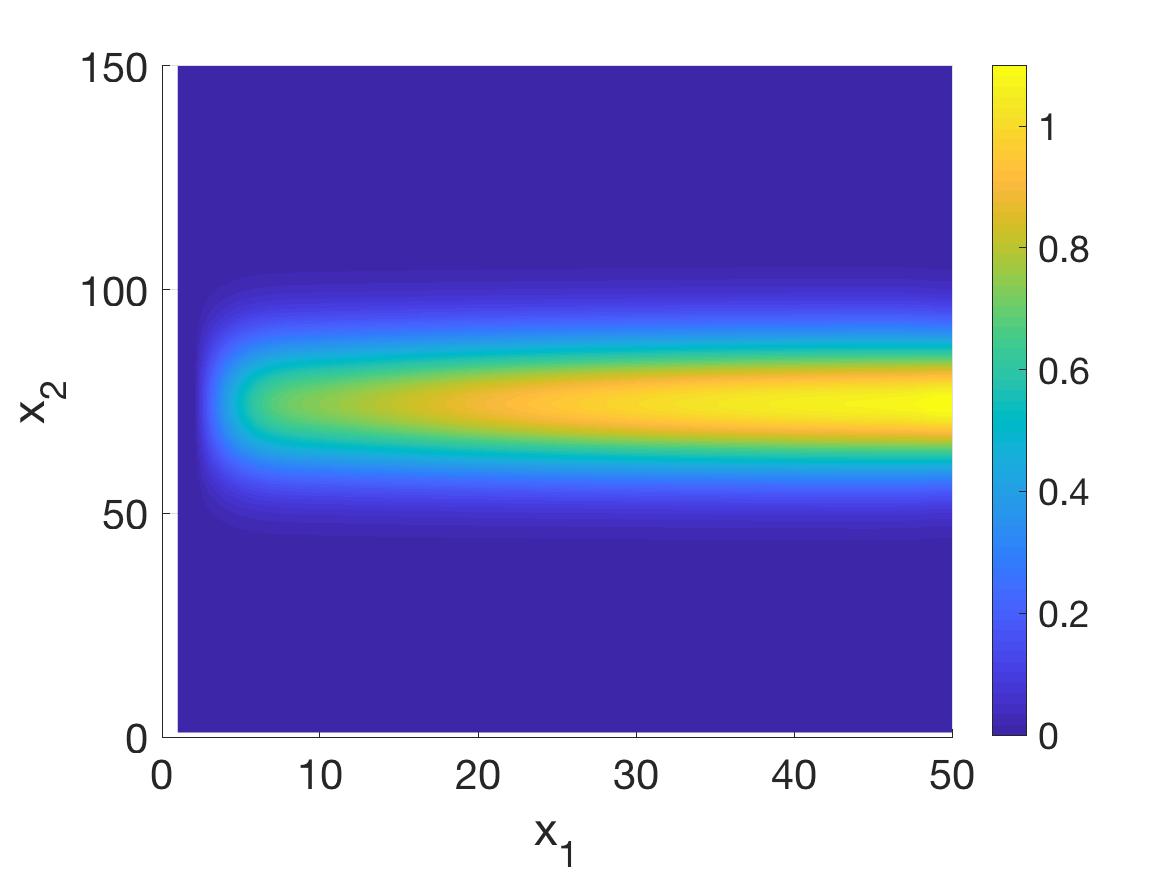} \hskip -5mm
\includegraphics[width=4cm]{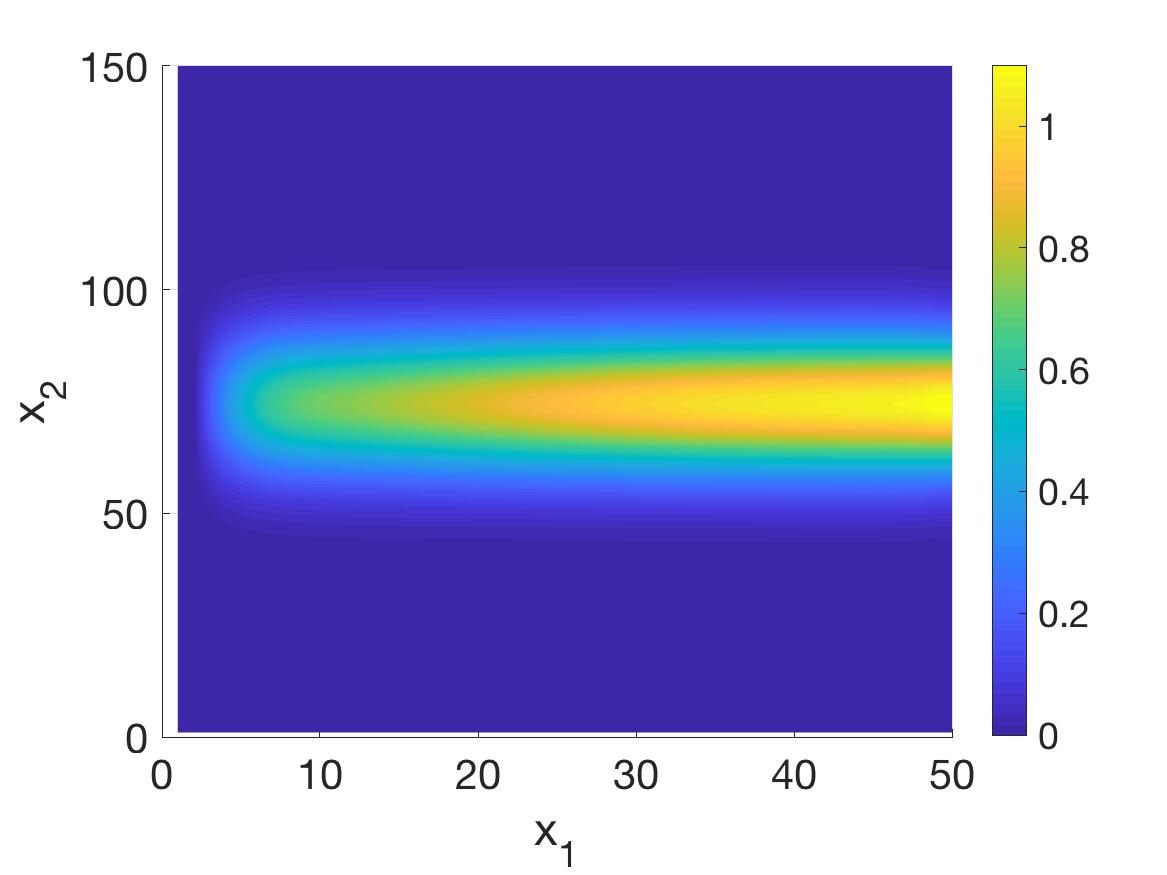} \\
(d) \hskip 3cm (e) \hskip 3cm (f) \\
\includegraphics[width=4cm]{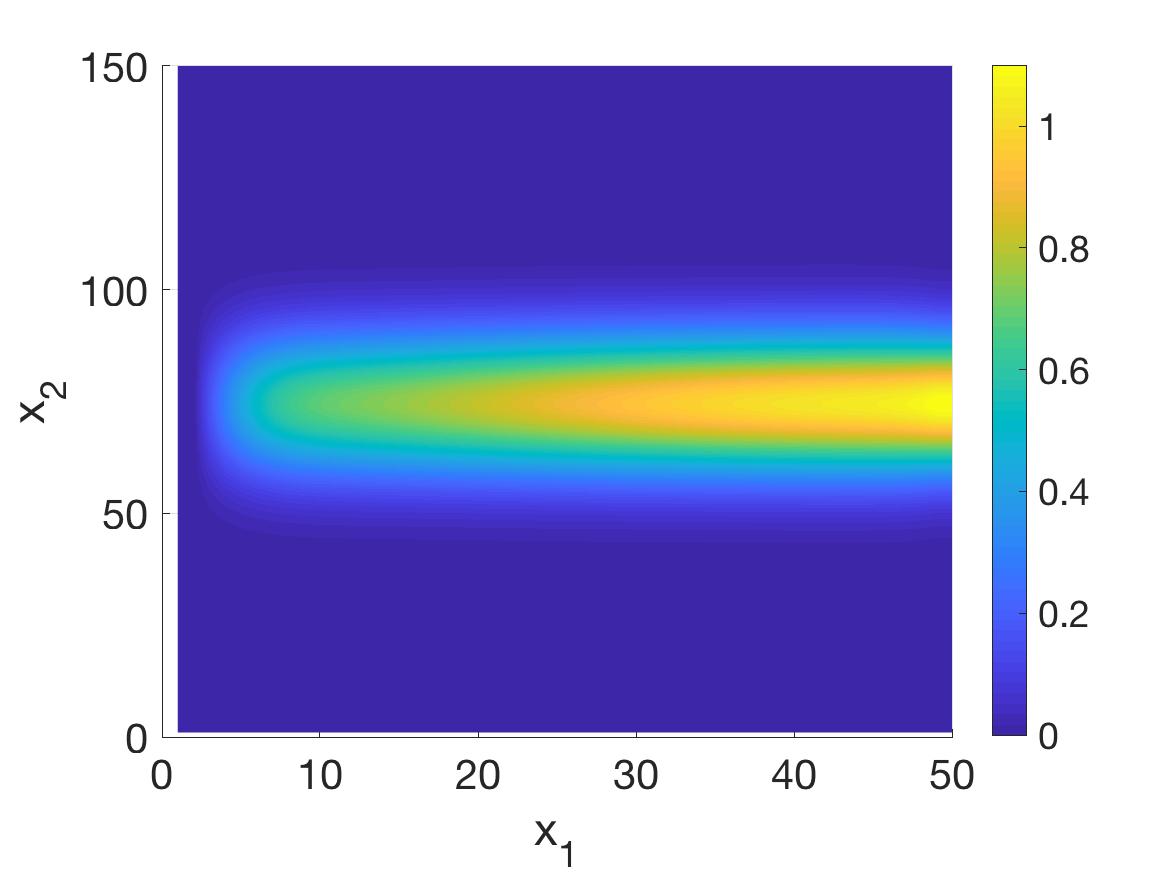} \hskip -5mm
\includegraphics[width=4cm]{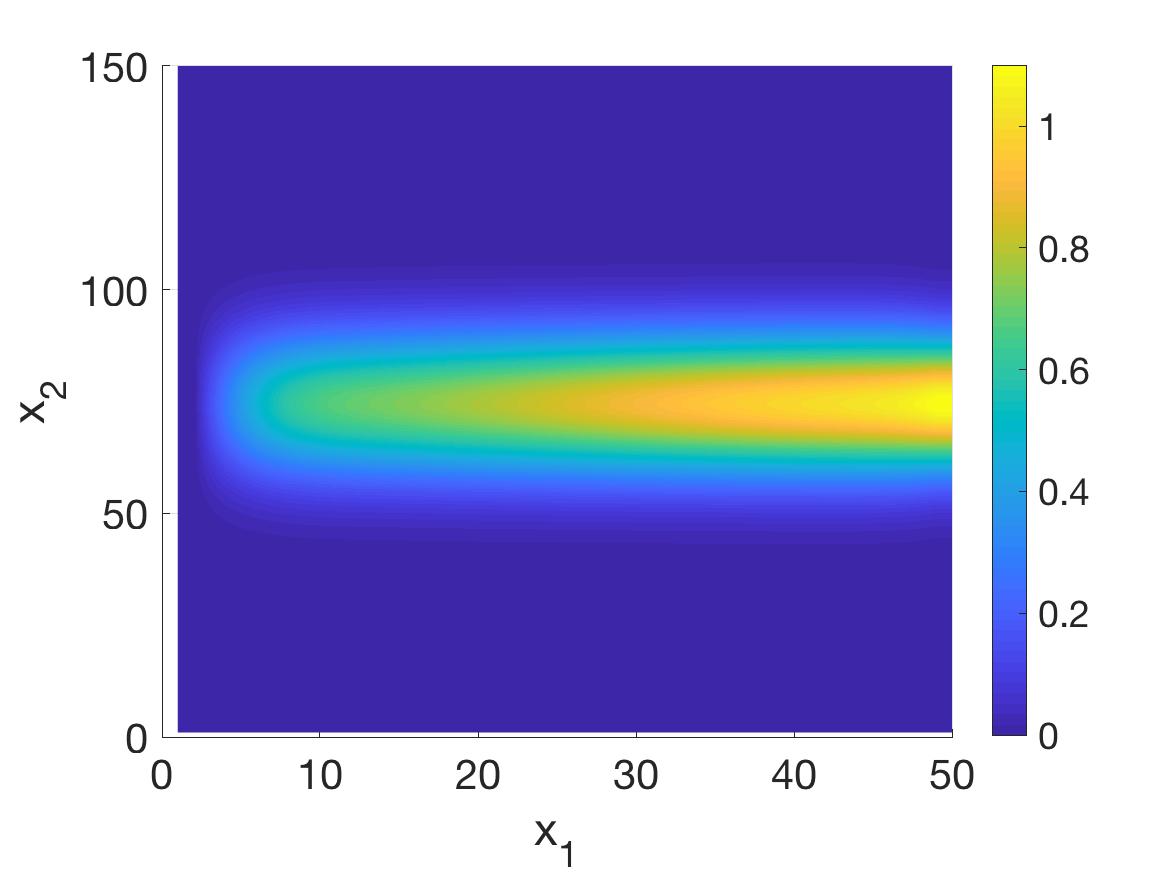} \hskip -5mm
\includegraphics[width=4cm]{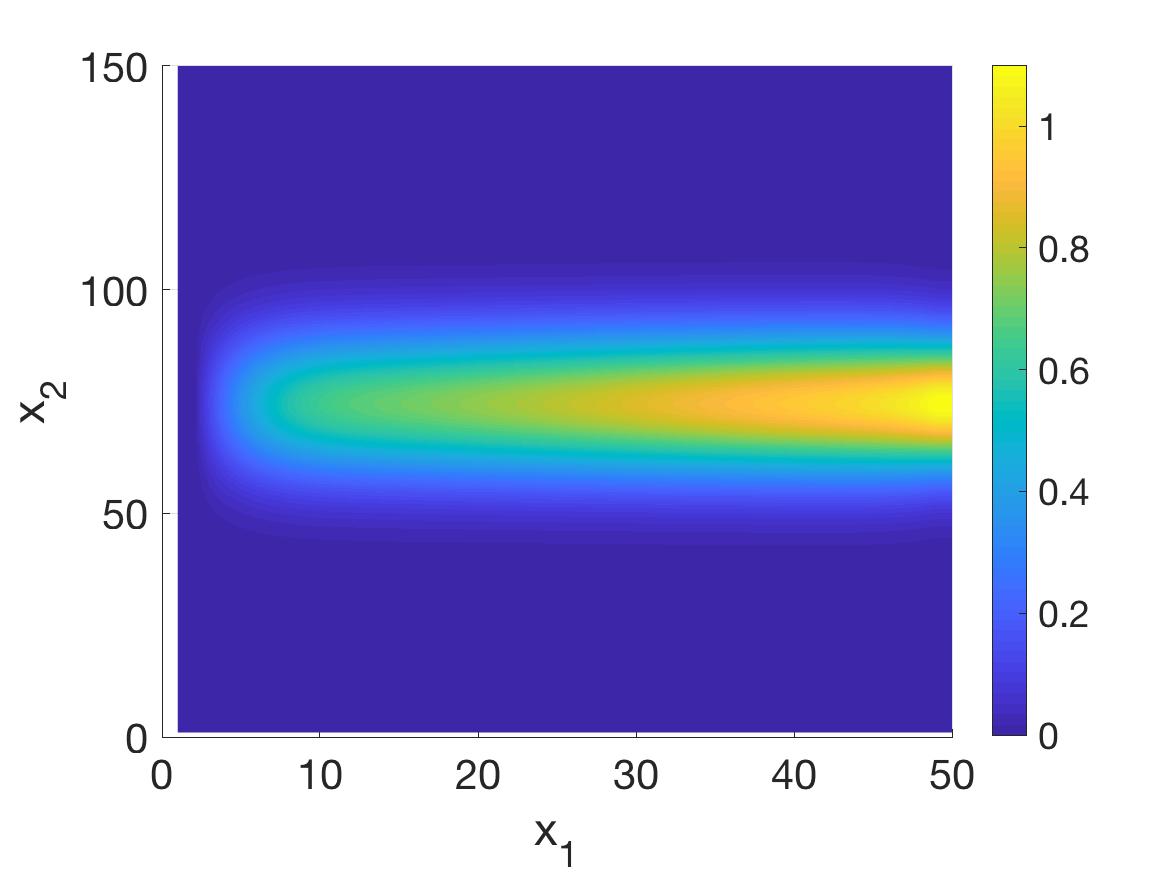} 
\caption{Concentration of tumor angiogenic factor obtained solving the two dimensional
reduced model at times
(a) $0.1145$, (b) $0.2290$, (c)  $0.3435$, (d) $0.4580$,
(e)  $0.5725$, (f) $0.6870$. 
The spread of the angiogenic factor shrinks as new vessel tips are formed and advance
 towards the hypoxic region.
}
\label{fig3}
\end{figure}

Figures \ref{fig2} and \ref{fig3} show the evolution of the marginal tip density
$\rho$ and the angiogenic factor concentration. Note that we only plot the
density of tips, not the full vessels already formed that lie behind them,
between the pre-existing vessel and the new tips.

\section{Conclusions}
\label{sec:conclusions}

We have introduced a high order positivity preserving scheme for asymptotic
reductions of kinetic models of angiogenesis, combining positivity preserving
high order WENO discretizations in space and SSP discretization in time.
The resulting schemes reproduce traveling patterns representing the formation
and spread of blood vessel tips in the direction of an hypoxic region, in need
of oxygen supply.
This approach provides building blocks to devise high order positivity 
preserving schemes for the full kinetic models, by resorting to fractional step 
splitting schemes \cite{split,split2}. The proper implementation of such schemes 
is at present a computational challenge, out of the scope of the present study.

\section{Appendix: Reconstruction from averages}
\label{sec:appendix}

For ease of the reader, we scketch here the WENO5 reconstruction procedure
of the required approximations of $u$ from the double averages 
$\overline{\overline{u}}_{ij},$ taken from the detailed description in \cite{shu}.
 This is done dimension by dimension. Let us first recall
how to reconstruct a function $u$ of one variable $x$ from the double averages 
$\overline{\overline{u}}_{i-2},$ $\overline{\overline{u}}_{i-1},$
$\overline{\overline{u}}_{i},$ $\overline{\overline{u}}_{i+1},$
$\overline{\overline{u}}_{i+2},$ 
in the variable $x$ at stencils ${\cal S}$ formed by intervals $I_{i-2},I_{i-1},I_{i},
I_{i+1},I_{i+2}$ with $I_i=[x_{i-{1\over 2}}, x_{i+{1\over 2}}]$. The coefficients of a polynomial
$p(x)= \sum_{\ell=0}^4 c_\ell \left({ x- x_i \over \delta x} \right)^\ell$
satisfying
\begin{eqnarray}
\overline{\overline{u}}_{j} = {1\over \delta x^2}
\int_{x_{j-{1\over 2}}}^{x_{j+{1\over 2}}} \int_{x-{\delta x\over 2}}^{x+{\delta x\over 2}}
p(\xi) \, d \xi dx, \quad  j=i-2,...,i+2, \label{weno1}
\end{eqnarray}
are given by
\begin{eqnarray} \begin{array}{l}
c_0=  {1\over 180} ( 2 \overline{\overline{u}}_{i-2}  - 23 \overline{\overline{u}}_{i-1} + 222 \overline{\overline{u}}_{i} -23 \overline{\overline{u}}_{i+1} 2 \overline{\overline{u}}_{i+2} ),
\\[1ex]
c_1= {1\over 8} (  \overline{\overline{u}}_{i-2}  - 6 \overline{\overline{u}}_{i-1} + 6 \overline{\overline{u}}_{i+1} - \overline{\overline{u}}_{i+2} ), 
\\[1ex]
c_2= {1\over 12} ( - \overline{\overline{u}}_{i-2} +10 \overline{\overline{u}}_{i-1} -18 \overline{\overline{u}}_{i} +10 \overline{\overline{u}}_{i+1} - \overline{\overline{u}}_{i+2} ), 
\\[1ex]
c_3= {1\over 12} ( - \overline{\overline{u}}_{i-2}  +2 \overline{\overline{u}}_{i-1} -2 \overline{\overline{u}}_{i+1} + \overline{\overline{u}}_{i+2} ),
\\[1ex]
c_4= {1\over 24} ( \overline{\overline{u}}_{i-2}  - 4 \overline{\overline{u}}_{i-1} + 6 \overline{\overline{u}}_{i} -4 \overline{\overline{u}}_{i+1} + \overline{\overline{u}}_{i+2} ).
\end{array} \label{coefs1d}
\end{eqnarray}
In terms of linear weights $d_m(x)$, $m=0,1,2$, this polynomial is decomposed
as
\begin{eqnarray} \begin{array}{ll}
p(x)= & \sum_{m=0}^2 d_m(x) p_m(x), \\[1ex]
p_0(x)= & {1 \over 12}(- \overline{\overline{u}}_{i-2} +2 \overline{\overline{u}}_{i-1} + 11 \overline{\overline{u}}_{i}) 
+ {1\over 2} (\overline{\overline{u}}_{i-2} - 4 \overline{\overline{u}}_{i-1} + 3 \overline{\overline{u}}_{i}) \left({ x- x_i \over \delta x} \right) \\
& + {1 \over 2} (\overline{\overline{u}}_{i-2} -2 \overline{\overline{u}}_{i-1}  +
\overline{\overline{u}}_{i}) \left({ x- x_i \over \delta x} \right)^2, \\[1ex]
p_1(x)= & {1 \over 12}(- \overline{\overline{u}}_{i-1} +14 \overline{\overline{u}}_{i} - \overline{\overline{u}}_{i+1}) 
+ {1\over 2} (-\overline{\overline{u}}_{i-1} + \overline{\overline{u}}_{i+1}) 
\left({ x- x_i \over \delta x} \right) \\
& + {1 \over 2} (\overline{\overline{u}}_{i-1} -2 \overline{\overline{u}}_{i}  +
\overline{\overline{u}}_{i+1}) \left({ x- x_i \over \delta x} \right)^2, \\[1ex]
p_2(x)= & {1 \over 12}(11\overline{\overline{u}}_{i} + 2 \overline{\overline{u}}_{i+1} - \overline{\overline{u}}_{i+2}) 
+ {1\over 2} (- 3\overline{\overline{u}}_{i} + 4  \overline{\overline{u}}_{i+1} -
\overline{\overline{u}}_{i+2})  \left({ x- x_i \over \delta x} \right) \\
& + {1 \over 2} (\overline{\overline{u}}_{i} -2 \overline{\overline{u}}_{i+1}  +
\overline{\overline{u}}_{i+2}) \left({ x- x_i \over \delta x} \right)^2.\\[1ex]
\end{array} \label{linearw} \end{eqnarray}
Imposing $\sum_{m=0}^2 d_m(x)=1$, these weights are listed in Tables
2.1 and 2.2 of Ref. \cite{shu} for the points $x_{i \pm {1\over 2}}^\alpha$,
$\alpha=1,2,3$ and $\tilde x_i^\alpha$, $\alpha=1,\ldots,5$ we use (see Table 
\ref{Table2} here).
Some of these linear weights are negative and may lead to instability. 
In practice, they are replaced by normalized nonlinear weights $\omega_m$. 
When at a point $x$, the three linear weights $d_m(x)$ are positive,
we set
\begin{eqnarray}\begin{array}{l}
\omega_m(x)= {\tilde w_m(x) \over \sum_{\ell=0}^2 \tilde w_\ell(x)},
\quad
\tilde w_m(x) = {d_m(x) \over (\beta_m + \varepsilon)^2}, \\[2ex]
\beta_0= {13 \over 12} (\overline{\overline{u}}_{i-2} - 2 \overline{\overline{u}}_{i-1}
+ \overline{\overline{u}}_{i})^2 + {1 \over 4} (\overline{\overline{u}}_{i-2}
 - 4 \overline{\overline{u}}_{i-1} + 3 \overline{\overline{u}}_{i})^2, \\[1ex]
\beta_1= {13 \over 12} (\overline{\overline{u}}_{i-1} - 2 \overline{\overline{u}}_{i}
+ \overline{\overline{u}}_{i+1})^2 + {1 \over 4} (\overline{\overline{u}}_{i-1}
 - 3 \overline{\overline{u}}_{i+1})^2, \\[1ex]
\beta_2 = {13 \over 12} (\overline{\overline{u}}_{i} - 2 \overline{\overline{u}}_{i+1}
+ \overline{\overline{u}}_{i+2})^2 + {1 \over 4} (3 \overline{\overline{u}}_{i}
 - 4 \overline{\overline{u}}_{i+1} + 3 \overline{\overline{u}}_{i+2})^2,
\end{array} \label{nonlinear} \end{eqnarray}
for $\varepsilon >0$ small (typically $\varepsilon =10^{-6}$)
and $m=0,1,2$. If any of the weights is negative at $x$, we define
$\omega_m(x)$ as follows
\begin{eqnarray}\begin{array}{l}
\omega_m(x) =\sigma^+(x) \omega_m^+(x) - \sigma^-(x) \omega_m^-(x), \quad
 \sigma^\pm(x) =\sum_{\ell=0}^2 \tilde  \gamma_\ell^{\pm}(x),
\\[1ex]
\omega_m^{\pm}(x)= {\tilde w_m^{\pm}(x) \over \sum_{\ell=0}^2 \tilde w_\ell^{\pm}(x)},
\quad
\tilde w_m^{\pm}(x) = {\gamma^{\pm}_m(x) \over (\beta_m + \varepsilon)^2}, 
\\[2ex]
\gamma_m^{\pm}(x) =  {\tilde \gamma_m^{\pm}(x) \over \sigma^\pm(x)}, \;
\tilde \gamma_m^+(x) = {1\over 2} (d_m(x)+ \theta |d_m(x)|), \;
\tilde \gamma_m^-(x) = \tilde \gamma_m^+(x) - d_m(x),
\end{array} \label{nonlinear2} \end{eqnarray}
for $\theta=3$ and $m=0,1,2$. Once this is done, the polynomial approximation
of $u$ in ${\cal S}$ obtained from the double averages is
\begin{eqnarray}
p(x)= \sum_{m=0}^2 \omega_m(x) d_m(x). \label{weno2}
\end{eqnarray}
The procedure to reconstruct a two dimensional function at the required points 
is as follows:
\begin{itemize}
\item For $u(x_{i+{1\over 2}}^\beta,\tilde{y}_j^\alpha)^-$ and 
$u(x_{i-{1\over 2}}^\beta,\tilde{y}_j^\alpha)^+$ 
\begin{itemize}
\item Fix $\tilde{y}_j^\alpha$, $\alpha=1,\ldots,5$, and apply (\ref{weno1})-(\ref{weno2}) 
to reconstruct one dimensional double averages 
$\overline{\overline{u}}_i(\tilde{y}_j^\alpha)$ in the $x$ direction from the
values $\overline{\overline{u}}_{i,j-2},$ $\overline{\overline{u}}_{i,j-1},$
$\overline{\overline{u}}_{i,j},$ $\overline{\overline{u}}_{i,j+1},$
$\overline{\overline{u}}_{i,j+2}.$ 
\item Next, apply  (\ref{weno1})-(\ref{weno2}) to reconstruct $u$ at the
points $(x_{i\pm{1\over 2}}^\beta,\tilde{y}_j^\alpha)$ for $\beta=1,2,3$
from the values values $\overline{\overline{u}}_{i-2}(\tilde{y}_j^\alpha),$ 
$\overline{\overline{u}}_{i-1} (\tilde{y}_j^\alpha),$ 
$\overline{\overline{u}}_{i}(\tilde{y}_j^\alpha),$ 
$\overline{\overline{u}}_{i+1}(\tilde{y}_j^\alpha),$ 
$\overline{\overline{u}}_{i+2}(\tilde{y}_j^\alpha).$ 
\end{itemize}
\item For $u(\tilde{x}_i^\alpha, y_{j+{1\over 2}}^\beta)^-$ and 
$u(\tilde{x}_i^\alpha, y_{j-{1\over 2}}^\beta)^+$ 
\begin{itemize}
\item Fix $\tilde{x}_i^\alpha$,$\alpha=1,\ldots,5$, and apply (\ref{weno1})-(\ref{weno2}) 
to reconstruct one dimensional double averages 
$\overline{\overline{u}}_j(\tilde{x}_i^\alpha)$ in the $y$ direction from the
values $\overline{\overline{u}}_{i-2,j},$ $\overline{\overline{u}}_{i-1,j},$
$\overline{\overline{u}}_{i,j},$ $\overline{\overline{u}}_{i+1,j},$
$\overline{\overline{u}}_{i+2,j}.$ 
\item Next, apply  (\ref{weno1})-(\ref{weno2}) to reconstruct $u$ at the
points $(\tilde{x}_i^\alpha, y_{j\pm {1\over 2}}^\beta)$ for $\beta=1,2,3$
from the values $\overline{\overline{u}}_{j-2}(\tilde{x}_i^\alpha),$ 
$\overline{\overline{u}}_{j-1}(\tilde{x}_i^\alpha),$
$\overline{\overline{u}}_{j}(\tilde{x}_i^\alpha),$ 
$\overline{\overline{u}}_{j+1}(\tilde{x}_i^\alpha),$
$\overline{\overline{u}}_{j+2}(\tilde{x}_i^\alpha).$  
\end{itemize}
\item For $u(x_{i},\tilde{y}_j^\alpha)$
\begin{itemize}
\item Fix ${x}_i$ and apply (\ref{weno1})-(\ref{weno2}) to
reconstruct one dimensional double averages 
$\overline{\overline{u}}_j(x_{i})$ in the $y$ direction from the
values $\overline{\overline{u}}_{i-2,j},$ $\overline{\overline{u}}_{i-1,j},$
$\overline{\overline{u}}_{i,j},$ $\overline{\overline{u}}_{i+1,j},$
$\overline{\overline{u}}_{i+2,j}.$ 
\item Next, apply  (\ref{weno1})-(\ref{weno2}) to reconstruct $u$ at the
points $(x_{i},\tilde{y}_j^\alpha)$ for $\alpha=1,\ldots,5$
from the values $\overline{\overline{u}}_{j-2}(x_{i}),$ 
$\overline{\overline{u}}_{j-1}(x_{i}),$ $\overline{\overline{u}}_{j}(x_{i}),$ 
$\overline{\overline{u}}_{j+1}(x_{i}),$ $\overline{\overline{u}}_{j+2}(x_{i}).$  
\end{itemize}
\item For $u(\tilde{x}_i^\alpha, y_{j})$
\begin{itemize}
\item Fix ${y}_j$ and apply (\ref{weno1})-(\ref{weno2}) to
reconstruct one dimensional double averages 
$\overline{\overline{u}}_i({y}_j)$ in the $x$ direction from the
values $\overline{\overline{u}}_{i,j-2},$ $\overline{\overline{u}}_{i,j-1},$
$\overline{\overline{u}}_{i,j},$ $\overline{\overline{u}}_{i,j+1},$
$\overline{\overline{u}}_{i,j+2}.$ 
\item Next, apply  (\ref{weno1})-(\ref{weno2}) to reconstruct $u$ at the
points $u(\tilde{x}_i^\alpha, y_{j})$ for $\alpha=1,\ldots,5$
from the values $\overline{\overline{u}}_{i-2}({y}_j),$ 
$\overline{\overline{u}}_{i-1}({y}_j),$ $\overline{\overline{u}}_{i}({y}_j),$ 
$\overline{\overline{u}}_{i+1}({y}_j),$ $\overline{\overline{u}}_{i+2}({y}_j).$  
\end{itemize}
\end{itemize}

 \begin{table}[ht]
\begin{center}\begin{tabular}{cccc}
 \hline
$x $ & $d_1(x)$ &$d_2(x)$& $d_3(x)$  \\  \hline
$x_{i-{1 \over 2}} - {\sqrt{15} \over 10} \delta x$ &${307+72 \sqrt{15} \over 960}$ & 
${8377-1542 \sqrt{15} \over 6720}$ & $ {173 (-11+3 \sqrt{15}) \over 3360}$ \\  
 \hline
$x_{i-{1 \over 2}}$ & ${341 \over 1200}$ &  ${337 \over 600}$  & ${37 \over 240}$ \\ 
 \hline
$x_{i-{1 \over 2}} + {\sqrt{15} \over 10} \delta x$  & ${307-72\sqrt{15} \over 960}$ & 
${8377+1542\sqrt{15} \over 6720}$ & $- {173 (11+3\sqrt{15}) \over 3360}$ \\  
 \hline
$x_{i+{1 \over 2}} - {\sqrt{15} \over 10} \delta x$ & $- {173 (11+3\sqrt{15}) \over 3360}$ & 
${8377+1542\sqrt{15} \over 6720}$ & ${307-72\sqrt{15} \over 960}$ \\ 
 \hline
$x_{i+{1 \over 2}}$  & ${37 \over 240}$ & ${337 \over 600}$ & ${341 \over 1200}$ \\  
 \hline
$x_{i+{1 \over 2}} + {\sqrt{15} \over 10} \delta x$  & ${173 (-11+3\sqrt{15}) \over 3360}$ & 
${8377-1542\sqrt{15} \over 6720}$ & ${307+72\sqrt{15} \over 960}$ \\ 
 \hline
$x_i - {\sqrt{10} \over 5} \delta x$ &$ {427+87\sqrt{15} \over 1590}$ & ${368\over 795}$ & 
${427-87\sqrt{15} \over 590}$ \\   
 \hline
$x_i - {\sqrt{10} \over 10} \delta x$ & ${29147-246\sqrt{15} \over 129360}$ & ${35533 \over 64680}$ & ${29147+246\sqrt{15} \over 129360}$ \\  
 \hline
$x_i$ & $- {2\over 15}$ &$ {19 \over 15}$ & $- {2 \over 15}$ \\  
 \hline
$x_i + {\sqrt{10} \over 10} \delta x$ & ${29147+246\sqrt{15} \over 129360}$ & ${35533 \over 64680}$ & ${29147-246\sqrt{15} \over 129360}$ \\  
 \hline
$x_i + {\sqrt{10} \over 5} \delta x$ & ${427-87\sqrt{15} \over 1590}$ & ${368 \over 795}$ & 
${427+87\sqrt{15} \over 1590} $ \\
 \hline
\end{tabular}
\end{center}
\caption{Linear weights $d_m(x)$ for typical point choices.} 
\label{table2}
\end{table}

\vskip 2mm

{\bf Acknowledgements.} 
This research has been partially supported by the FEDER/ Ministerio de Ciencia, 
Innovación y Universidades -Agencia Estatal de Investigación grant No. 
MTM2017-84446-C2-1-R.

\end{document}